\documentclass{tac}

\usepackage{amssymb}
\usepackage{amsmath}
\usepackage{tikz, tikz-cd}

\usepackage{diagrams}

\newtheorem{theorem}{Theorem}[section]

\newtheorem{definition}[theorem]{Definition}

\newtheorem{proposition}[theorem]{Proposition}

\title{A classifying localic category for locally compact locales with application to the Axiom of Infinity}

\author{Dr Christopher F. Townsend}
\copyrightyear{2024}

\begin{document}
\maketitle

  \begin{abstract}

    For an internal category $\mathbb{C}$ in a cartesian category $\mathcal{C}$ we define, naturally in objects $X$ of $\mathcal{C}$, $Prin_{\mathbb{C}}(X)$. This is a category whose objects are principal $c \mathbb{C}$-bundles over $X$ and whose morphisms are principal $c(\mathbb{C}^{\uparrow})$-bundles. Here $c(\_)$ denotes taking the core groupoid of a category (same objects but only isomorphisms as morphisms) and $\mathbb{C}^{\uparrow}$ is the arrow category of $\mathbb{C}$ (objects morphisms, morphisms commuting squares). We show that $X \mapsto Prin_{\mathbb{C}}(X)$ is a stack of categories and call stacks of this sort lax-geometric. We then provide two sufficient conditions for a stack to be lax-geometric and use them to prove that the pseudo-functor $X \mapsto \mathbf{LK}_{Sh(X)}$ on the category of locales $\mathbf{Loc}$ is a lax-geometric stack. Here $\mathbf{LK}_{Sh(X)}$ is the category of locally compact locales in the topos of sheaves over $X$, $Sh(X)$. Therefore there exists a localic category $\mathbb{C}_{\mathbf{LK}}$ such that $\mathbf{LK}_{Sh(X)} \simeq Prin_{\mathbb{C}_{\mathbf{LK}}}(X)$ naturally for every locale $X$. 

    We then show how this can be used to give a new localic characterisation of the Axiom of Infinity. 

  \end{abstract}

\bigskip

{\bf This is a copy of the poster presented to CT2024. Proofs are outline only.}

  \section{A problem with morphisms between Principal Bundles}

A problem with principal bundles is that all morphisms between them are isomorphisms. So if we are hoping to use principal bundles to represent stacks of categories rather just stacks of groupoids we are stuck. Here we are proposing a way round this by defining a category of principal bundles relative to any internal category. We then ignore the generality (from groupoids to categories) at the level of objects but use it to define morphisms. The idea is encapsulated in the following definition.

\begin{definition}
For an internal category $\mathbb{C}$ in a cartesian category $\mathcal{C}$, define for each object $X$ of $\mathcal{C}$, the data $ Prin_{\mathbb{C}}(X)$ by: 

1. {\bf Objects}. The objects are principal $c\mathbb{C}$-bundles over $X$.

2. {\bf Morphisms}. For any two principal $c\mathbb{C}$-bundles $P_0$ and $P_1$ over $X$, a morphism from $P_0$ to $P_1$ is a principal $c(\mathbb{C}^{\uparrow})$-bundle $P$ over $X$ together with two principal $c\mathbb{C}$-bundle maps over $X$, $\alpha_0: P_0 \rTo \Sigma_{D_0} P$ and $\alpha_1: \Sigma_{D_1} P \rTo P_1$. Here $\Sigma_{D_i}P$ is obtained from the cocycle of $P$ by composition with the domain/codomain internal functors $D_i:c(\mathbb{C})^{\uparrow} \rTo c\mathbb{C}$. 

3. {\bf Composition}. The composition $P' \circ P$ of $P_0 \rTo^P P_1$ and $P_1 \rTo^{P'} P_2$ is the principal bundle associated with the cocycle
\begin{eqnarray*}
\psi^{P' \circ P}:\mathbb{P} \times_{\mathbb{X}} \mathbb{P}_1 \times_{\mathbb{X}} \mathbb{P'} \rTo {(\psi^P \pi_1, \gamma, \psi^{P'} \pi_3)}  c(\mathbb{C})^{\uparrow}\times_{c\mathbb{C}} c(\mathbb{C}^{\uparrow}) \rTo^{M} c(\mathbb{C}^{\uparrow})
\end{eqnarray*}
where $\gamma$ is the internal natural transformation $(\alpha_1)_{\pi_{12}}$ followed by $(\alpha'_0)^{-1}_{\pi_{23}}$. 

4. {\bf 2-cells}. A 2-cell $P \rTo^\rho \tilde{P}$ consists of a $c(\mathbb{C}^{\uparrow})$-bundle morphism over $X$ compatible with the isomorphisms defining $P$ and $\tilde{P}$; i.e. $\tilde{\alpha}_0 \Sigma_{D_0}\rho = \alpha_0$ and $\tilde{\alpha}_1 \Sigma_{D_1}\rho = \alpha_1$.

\end{definition}
\section{We have a well defined category}
\begin{proposition}
    $ Prin_{\mathbb{C}}(X)$ is a bicategory, with all 2-cells isomorphisms and there is at most one 2-cell between any two morphisms.  
    \end{proposition}

Because there is at most one isomorphism between morphisms we can treat $Prin_{\mathbb{C}}(X)$ as a category.

\section{We haven't `thrown the baby out with the bathwater'}
Crucially we haven't introduced or lost any of the usual morphisms (isomorphisms) between principal bundles because by taking the core of $Prin_{\mathbb{C}}(X)$ we get, up to equivalence, the usual category of principal bundles over the core of $\mathbb{C}$: 
\begin{proposition}
Given an internal category $\mathbb{C}$ in a cartesian category $\mathcal{C}$,
\begin{eqnarray*}
    Prin_{c \mathbb{C}}(X) \simeq c(Prin_{\mathbb{C}}(X))
\end{eqnarray*}
naturally in objects $X$ of $\mathcal{C}$.
\end{proposition}

\section{Introducing `lax' geometric stacks}
Now that we have defined new morphisms between principal bundles, it is natural to isolate the categories $Prin_{\mathbb{C}}(X)$ via a stack.

\begin{definition}
    A pseudo-functor $M: \mathcal{C}^{op} \rTo \frak{CAT}$ is a \emph{lax-geometric stack} if it is naturally equivalent to $Prin_{\mathbb{C}}(\_)$ for some category $\mathbb{C}$ internal to $\mathcal{C}$.
 
\end{definition}

The following is a lax version of a well established characterisation of geometric stacks: 
\begin{proposition}\label{III}{\bf A simple characterisation of lax-geometric stacks.}
    A pseudo-functor $M: \mathcal{C}^{op} \rTo \frak{CAT}$ is a \emph{lax-geometric stack} if and only if it is a stack and:  

(i) There exists an object $C$ in $M(C_0)$ for some object $C_0$ of $\mathcal{C}$ such that for any other object $A$ of $M(X)$, there exists a cover $p_A:Y \rTo X$, a morphism $q_A:Y \rTo C_0$ and an isomorphism $M(p_A)(A) \cong M(q_A)(C_0)$.

(ii) There exists an object $C_1$ and a bijection
\begin{eqnarray*}
\mathcal{C}(X,C_1) \cong \{(f,\theta,g) | f,g:X \rTo G_0, \theta:M(f)(C) \rTo M(g)(C) \} 
\end{eqnarray*}
natural in $X$.
\end{proposition}

Our {\bf Main Theorem} is that these conditions hold for the stack of locally compact locales:
\begin{theorem}
    The stack $\mathfrak{LK}$ defined by 
    \begin{eqnarray*}
        \mathbf{Loc}^{op} & \rTo & \mathfrak{CAT}\\
        X & \mapsto & \mathbf{LK}_{Sh(X)}
    \end{eqnarray*}
    is lax geometric. 
\end{theorem}

 \begin{proof}
    Check (i) and (ii) of Proposition \ref{III}. Checking (i) follows the now well established trick of seeing a particular type of frame (here a continuous frame) as always coming from the completion of a model of a geometric theory (see for example Proposition 4.4 of \cite{HenryTowClass}). The relevant geometric theory to use is an information system (\cite{vicinfosys}) augmented with the notion of strength taken from \cite{Kaw}; the completion is taking rounded ideals.
    
    For (ii), note that the $C$ derived in the verification of (i) is locally compact over some locale $C_0$. Therefore the localic pullback $\pi_i^*C$ is locally compact where $\pi_i: C_0 \times C_0 \rTo C_0$ for $i = 1,2$. But locally compact locales are exponentiable. Define $C_1 = (\pi_2^*C)^{\pi_1^*C}$; the exponentiation is in the category of locales over $C_0 \times C_0$.    
    \end{proof}
    
\section{How is this related to the Axiom of Infinity?}

Andreas Blass proved (\cite{blass}) that an elementary topos $\mathcal{E}$ has a natural numbers object if and only if there is an object classifier relative to $\mathcal{E}$ (see Theorem 4.2.11 of \cite{Elephant} for a textbook account). 

Having a natural numbers object is, in the context of topos theory, essentially the same thing as satisfying the Axiom of Infinity. So Blass's remarkable result effectively provides a novel characterisation of the Axiom of Infinity. 

But having an object classifier is really the same thing as knowing that $X \mapsto Sh(X)$ is a geometric stack. This can be seen by \cite{BDesc}, which is applicable since by Joyal and Tierney, \cite{JoyT}, we know that the object classifier must be of the form $B\mathbb{G}$ for some localic groupoid $\mathbb{G}$. With this viewpoint what we can derive from our Main Theorem is: 

\begin{theorem}{\bf Main Theorem as Axiom of Infinity.}
For any elementary topos $\mathcal{E}$, $\mathcal{E}$ has a natural numbers object if and only if the pseudo-functor
\begin{eqnarray*}
    \mathfrak{LK}:\mathbf{Loc}_{\mathcal{E}} & \rTo & \mathfrak{CAT}\\
    X & \mapsto & \mathbf{LK}_{Sh_{\mathcal{E}}(X)}\\
\end{eqnarray*}
is a lax geometric stack with a localic category $\mathbb{C}_{\mathbf{LK}}$ whose core is stably Frobenius.
\end{theorem}
We say that an internal groupoid $\mathbb{G}$ is \emph{stably Frobenius} if it has a stably Frobenius connected components adjunction. For any cartesian $\mathcal{C}$ there is a 2-categorical functor $SFGpd_{\mathcal{C}} \rTo T_{\mathcal{C}}$ that universally inverts essential equivalences. Here $T_{\mathcal{C}}$ is the 2-category whose objects are stably Frobenius connected component adjunctions and whose morphisms are stably Frobenius adjunctions over $\mathcal{C}$. The paper \cite{TowHS} provides a lot of technical details about the 2-category $T_{\mathcal{C}}$; the key properties for our context being that all 2-cells are isomorphisms and it has finite 2-limits. We use this information to give an outline proof of our characterisation of the Axiom of Infinity:   

\begin{proof}{\emph{Outline Proof of the relationship to the Axiom of Infinity}.}

    The main theorem shows how $\mathbb{C}_{\mathbf{LK}}$ is constructed. Further it is the completion of an open localic groupoid that classifies a geometric theory. So we can use the fact the any open localic groupoid is stably Frobenius to show that $c\mathbb{C}_{\mathbf{LK}}$ is stably Frobenius. 
    
    For the other direction note that the stack $c\mathfrak{LK}$ can be extended to a pseudo functor $T_{\mathbf{Loc}_{\mathcal{E}}}^{op} \rTo \mathfrak{GPD}$ which is then representable by assumption. But moreover $c\mathfrak{LK}^{\uparrow}$ can also be so extended and is also representable; it is represented by the core of the arrow category of $\mathbb{C}_{\mathbf{LK}}$. Now there are two natural transformations $\mathbb{S}^{(\_)},P_L(\_):c\mathfrak{LK} \rTo c\mathfrak{LK}$ which therefore correspond to two maps in $T_{\mathbf{Loc}_{\mathcal{E}}}$; their equalizer must represent discrete locales because a locale $X$ is discrete if and only if $\mathbb{S}^X \cong P_L(X)$. Because we have classified discrete locales and we have a classifying localic groupoid for morphisms between locales we can construct a classifying groupoid for morphisms between discrete locales (because $T_{\mathcal{C}}$ has limits). Moreover, we can construct a classifying groupoid for objects of the form $1 \rTo A \rTo A$ where $A$ is discrete. But then, even though we don't know that the groupoid constructed is \'{e}tale complete, it is still possible to construct a direct image from its connected components adjunction. This is done by exploiting how natural transformations between frames correspond to dcpo homomorphisms, \cite{victow}. Then the argument of Theorem 4.2.11 of \cite{Elephant} can be repeated to prove that $\mathcal{E}$ has a natural numbers object.   
\end{proof}

The above represents a proposed Axiom of Infinity for locales; that is, it could be explored in the context of axiomatic approaches to locale theory (e.g. \cite{TowIdemp}). It can be used to show aspects of Hyland's result on `locally compact = exponentiable' even without any underlying set theory. It is symmetric under open/proper duality and implies the existence of a classifying category for both discrete objects and compact Hausdorff objects (c.f. \cite{HenryTowClass}). We conjecture that because of this symmetry it cannot imply that a discrete axiom of infinity holds relative to axiomatic locale theory.

\end{document}